\documentclass[12pt]{amsart}
\usepackage{amsfonts,amssymb,amscd,amsmath}
\usepackage{latexsym}
\usepackage{tocvsec2}
\usepackage{mathrsfs}
\usepackage[all]{xy}

\usepackage[
hyperindex=true,pagebackref=true,bookmarks=true,
colorlinks=true,linkcolor=blue,citecolor=red]
{hyperref}
\def\~{{\rm --}}

\settocdepth{subsubsection}

\renewcommand{\tilde}{\widetilde}
\renewcommand{\hat}{\widehat}

\newcommand{\Z}{{\mathbb Z}}
\newcommand{\Q}{{\mathbb Q}}

\newcommand{\R}{{\mathbb R}}

\def\HH{\mbox{${\mathcal H}$\kern-5.2pt${\mathcal H}$}}

\newtheorem{theorem}{Theorem}[section]

\newtheorem{theorem }{Theorem}[section]
\newtheorem{maintheorem }[theorem]{Main Theorem}
\newtheorem{proposition }[theorem]{Proposition}
\newtheorem{definition }[theorem]{Definition}
\newtheorem{lemma }[theorem]{Lemma}
\newtheorem{corollary }[theorem]{Corollary}
\newtheorem{notation }[theorem]{Notation}
\newtheorem{remark }[theorem]{Remark}
\newtheorem{example }[theorem]{Example}

\newtheorem{ maintheorem }[theorem]{Main Theorem}
\newtheorem{ theorem}{Theorem}[section]
\newtheorem{ proposition}[theorem]{Proposition}
\newtheorem{ definition}[theorem]{Definition}
\newtheorem{ lemma}[theorem]{Lemma}
\newtheorem{ corollary}[theorem]{Corollary}
\newtheorem{ notation}[theorem]{Notation}
\newtheorem{ remark}[theorem]{Remark}
\newtheorem{ example}[theorem]{Example}

\def\for{\  \hbox{ for } \ }

\def\where{\  \hbox{ where } \ }
\def\and{\  \hbox{ and } \ }
\def\and{\  \hbox{ and } \ }

\def\equal{\stackrel{\,\mathbf{def}}{= \kern-3pt =}}

\def\la{\lambda}
\def\La{\Lambda}
\def\om{\omega}

\def\th{\theta}
\def\al{\alpha}
\def\be{\beta}

\def\de{\delta}

\def\si{\sigma}

\def\Ga{\Gamma}
\def\ze{\zeta}

\def\vep{\varepsilon}

\def\tal{\tilde{\alpha}}

\def\tR{\tilde R}

\def\hw{\widehat{w}}
\def\hW{\widehat{W}}

\def\0{\mathbf{0}}

\def\çF{\mathcal{F}}

\def\n{\mathcal{N}}

\def\p{\mathcal{P}}

\def\c{\mathcal{C}}

\def\v{\mathcal{V}}

\def\lan{\langle}

\def\ran{\rangle}

\def\lng{\hbox{\rm{\tiny lng}}}
\def\sht{\hbox{\rm{\tiny sht}}}

\newcommand{\sq}{\phantom{1}\hfill$\qed$}

\newcommand{\lr}{\langle}
\newcommand{\rr}{\rangle}

\def\HH{\mathfrak{H}}

\def\HH{\hbox{${\mathcal H}$\kern-5.2pt${\mathcal H}$}}
\def\HHH{\hbox{${\mathbb H}$\kern-4.2pt${\mathbb H}$}}

\font\smm=msbm10 at 12pt 
\def\symbol#1{\hbox{\smm #1}}
\def\lsmash{{\symbol n}}

\def\#{\sharp}

\begin{document}
\newcommand{\comment}[1]{}

\begin{abstract}
This paper is devoted to the Harish-Chandra-type 
decomposition of the global
nonsymmetric spherical functions in terms
of their asymptotic expansions and the $q,t$\~generalization
of the celebrated $c$\~function. This is for any
reduced root systems in the $q,t$\~setting; we 
pay special attention to the case of $A_1$, 
where this decomposition is very explicit.
\end{abstract}

\title[Harish-Chandra theory 
of global functions]
{On Harish-Chandra theory 
of global nonsymmetric functions}
\author[Ivan Cherednik]{Ivan Cherednik $^\dag$}

\address[I. Cherednik]{Department of Mathematics, UNC
Chapel Hill, North Carolina 27599, USA\\
chered@email.unc.edu}

\def\bysame{{\bf --- }}
\def\~{{\bf --}}
\renewcommand{\tilde}{\widetilde}
\renewcommand{\hat}{\widehat}
\newcommand{\dagx}{\hbox{\tiny\mathversion{bold}$\dag$}}
\newcommand{\ddagx}{\hbox{\tiny\mathversion{bold}$\ddag$}}

\thanks{$^\dag$ \today.
\ \ \ Partially supported by NSF grant
DMS--1363138 and the Simons Foundation.}
\maketitle

\smallskip

{\small
\centerline{{\bf MSC} (2010): 20C08, 33C52, 33C67, 33D52, 33D67}
} 

\renewcommand{\baselinestretch}{1.2}
{
\smallskip
\tableofcontents
\smallskip
}
\renewcommand{\baselinestretch}{1.2}
\vfill\eject

\renewcommand{\natural}{\wr}
\renewcommand{\baselinestretch}{1.2}

\setcounter{section}{-1}
\setcounter{equation}{0}
\section{\sc Introduction}
In this paper, we obtain the Harish-Chandra-type 
asymptotic decomposition of the {\em global nonsymmetric 
spherical\,} $q,t$\~function $G(X,\La)$ from \cite{C5}. 
Given a Weyl chamber $\,\c\,$,
$G(X,\La)$ is represented as a weighted sum of its
asymptotic expansions for $\,\c \ni\Re(x)\sim
\infty\,$ for $\La$ and its translation by
$w\in W$ for the Weyl group $W$, where we set $X=q^x$. 
The weight functions are given in 
terms of  $\si(\La)$ from \cite{ChW},
a $q,t$\~generalization of the celebrated 
Harish-Chandra $c$\~function. Only reduced root 
systems will be considered. 

The role of $\si(\La)$ here is similar to that in the symmetric case; 
see \cite{ChW,Sto2} for the corresponding Harish-Chandra 
theory. However the nonsymmetric asymptotic
decomposition is not a $W$\~symmetrization in the $\La$\~space
since $G(X,\La)$ is not $W$\~invariant. The asymptotic 
expansions of $G(X,w(\La))$ are not connected with each 
other for $w\in W$ in any direct way. They form a 
{\em $W$\~spinor\,} solving the corresponding difference 
Dunkl eigenvalue problem in the terminology of 
\cite{C102,ChO1,ChO2}. This is parallel to the 
(differential) nonsymmetric theory in \cite{Op}. 

Using the (spinor) Dunkl eigenvalue problem is of obvious
significance here. It provides the asymptotic expansions and 
the existence of the required decomposition with 
undetermined coefficients. The first of these applications 
can be actually replaced by using the technique of intertwining 
operators from \cite{C1,HHL,RY,OS}, the second is very classical. 
The Dunkl eigenvalue problem (and $W$\~spinors) will be discussed in 
our further paper(s); we do not need it too much in this
particular paper, as well as the related DAHA theory.
 
We try to make this work short and as focused as 
possible. Also, to simplify the setup, only the negative 
Weyl chamber $\,\c_-\ni \Re(x)\,$ will be considered, where
$\,\Re (x,\al)<0$ for all positive roots $\al$. Actually it
suffices to establish the $\si$\~decomposition only in some
open $X$\~set (all involved function are meromorphic).
Arbitrary Weyl chambers will be hopefully considered in our further
works.
\smallskip

The function $G(X,\La)$ is given by a series that converges 
anywhere. However its asymptotic expansion in
$\,\c_-\ni \Re(x)\sim \infty\,$ is only meromorphic 
and has certain finite radius of convergence 
(dependent on $q,t$). Here $\La$
can be arbitrary avoiding singularities. For sufficiently
large negative $\Re(x)$, the cancelation of the 
$\La$\~singularities in the $\si$\~decomposition of 
$G(X,\La)$ is an impressive application of the theory of global 
functions. We note that generally 
there are no a priori ways for establishing such a 
cancelation in the (differential) Heckman-Opdam theory of the 
hypergeometric function \cite{HO,Op}. Now it can be (potentially)
deduced from the difference theory by taking the limit $q\to 1$.

We essentially follow \cite{ChW}, switching from the 
asymptotic expansions of the symmetric Macdonald polynomials 
there to those for the nonsymmetric Macdonald polynomials. 
The method from \cite{Sto2} is restricted to the symmetric theory; 
the difference Cherednik-Matsuo isomorphism theorem is used there, 
connecting the eigenvalue problem for 
the Macdonald operators  with the difference 
{\em AKZ-system\,} \cite{C101,C102,MS}.
This approach is not needed in the 
nonsymmetric theory because the Dunkl operators can be used 
instead. As a matter of fact, the (spinor) Dunkl eigenvalue problem
is the key in the justification of the Cherednik-Matsuo 
theorems. Note that \cite{Sto2} includes non-reduced
root systems (the present one is restricted to the reduced case).
\medskip

At almost any levels, the nonsymmetric direction simplifies,
clarifies and generalizes the symmetric one. This is 
especially true in the difference theory, which already 
significantly changed the classical harmonic analysis on 
symmetric spaces. This is fully applicable to this paper. However we
do not have any geometric interpretation of the nonsymmetric 
$\si$\~decomposition formula at the moment. Moreover, the 
geometric meaning of the $E$\~polynomials themselves for generic 
$q,t$ remains essentially unknown; the symmetrization is generally
needed here (apart from several limiting cases).
   
We pay special attention to $A_1$ in this paper, 
where our formula extends that from 
\cite{ChO1}. Here the decomposition is very explicit and the
connection with the classical {\em basic hypergeometric 
function\,} can be readily seen. Also, such an explicit
formula makes it possible to analyze its behavior at $|q|=1$.
This was touched upon in \cite{ChO1} ($A_1$, the symmetric case),
but will not be discussed in the present paper. 

\smallskip
{\bf Acknowledgements.}
The author is thankful to M.~Duflo, D.~Kazhdan, 
E.~Opdam, J.~Stokman  for various discussions of the 
Harish-Chandra theory and M.~Finkelberg for a discussion
of \cite{BFS}. This paper was stimulated by author's talks
at RIMS and in Tokyo University in 2014, as well as
B.~Feigin's question concerning the existence of $W$\~orbit-sum 
formulas for the nonsymmetric Macdonald polynomials 
(see (\ref{SphPsixx}) at the end). 
The paper was mainly written at RIMS; the author thanks Hiraku 
Nakajima and RIMS for the invitation and hospitality and
the Simons Foundation (which made this visit possible).
\medskip

\setcounter{equation}{0}
\section{\sc The case of \texorpdfstring{$A_1$}{A1}}
We will begin this paper with the case of $A_1$, when the
formulas for the $E$\~polynomials and the related
functions are (exceptionally) explicit. 
Let us recall the definition of the global function
from \cite{ChW} in the case of $A_1$. See Theorem 5.4 there
and also \cite{ChO1}.

The constant term functional, the coefficient of $X^0$
of a Laurent series $f$ or a polynomial in terms of $X^{\pm 1}$
will be denoted by $\lan\cdot\ran$. The $\mu$\~function,
the measure that makes the $E$\~polynomials orthogonal,
is the following truncated theta function:
\begin{align}\label{mucont}
&\mu(X;q,t)\equal\prod_{j=0}^\infty
\frac{(1-q^jX^2)(1-q^{j+1}X^{-2})}
{(1-tq^jX^2)(1-tq^{j+1}X^{-2})}.
\end{align}
We will mainly need its renormalization
\begin{align}\label{mucontone}
&\mu_\circ\equal\mu/\lr\mu\rr= 1+\frac{t-1}{1-qt}(X^2+qX^{-2})
+\ldots\ ,\\
&\where \lr \mu \rr\ =\ \prod_{j=1}^\infty \frac{(1-tq^j)^2}
{(1-t^2q^j)(1-q^j)}.\notag
\end{align}

The series $\mu_\circ$ is $\ast$\~invariant for the conjugation
$$
X^\ast=X^{-1},\ (q^{\,1/2})^\ast=q^{-1/2},\
(t^{\,1/2})^\ast=t^{-1/2}.
$$

The Demazure-Lusztig operator $T$ and the
difference Dunkl operator $Y$ are as follows:
\begin{align*}
&T\ =\ t^{1/2}s+\frac{t^{1/2}-t^{-1/2}}{X^{2}-1}
\cdot (s-1),\ \, Y\ =\ s\,\Ga\, T,\\
&s(X^n)=X^{-n},\ \Ga(X^n)=q^{n/2}X^{n} \for n\in\Z.
\end{align*}
They naturally act in the polynomial DAHA representation,
which is  $\v\equal\Z[q^{\pm1/2},t^{\pm1/2}][X^{\pm1}]$.
We will sometimes set $X=q^x$; then
$s(x)=-x,\, \Ga(f(x))=f(x+1/2)$.

The {\em nonsymmetric Macdonald polynomials},
also commonly called {\em $E$\~polynomials},
are uniquely determined  from the 
eigenvalue problem
\begin{align}\label{nonsymp}
Y(E_{n})=q^{-n_{\#}}E_{n}\for n\in \Z,\ t\equal q^k,&&&&\\
n_{\#}=\left\{\begin{array}{ccc}\frac{n+k}{2}
&\hbox{for}  & n>0, \\\frac{n-k}{2} &\hbox{for}  
& n\le 0 \end{array}\right\},
\text{\, note that }\,  0_{\#}=-\frac{k}{2},
\end{align}
where the normalization is
$E_{n}=X^{n}+\text{ ``lower terms'' }.$
By ``lower terms'', we mean
polynomials in terms of $X^{\pm m}$ as $|m|<n$
and, additionally,  $X^{|n|}$ for $n<0$.
It gives a filtration in $\v$
preserved by $Y$. These polynomials (for any reduced
root systems) are due to Opdam  in the differential setting 
(he mentions a contribution of Heckman) and
Macdonald for integral $k\ge 0$; see \cite{Op,M4,C4,C1}.

Obviously,  $E_{0}=1,\, E_{1}=X$.
Let us provide the formulas for the $E$\~polynomials ($n>0$):
\begin{align}\label{exacte-x}
E_{-n}=X^{-n}&+X^n\frac{1-t}{1-tq^n}+
\sum_{j=1}^{[n/2]}X^{2j-n}\,\prod_{i=0}^{j-1}
\frac{(1-q^{n-i})}{(1-q^{1+i})}\,\frac{(1-tq^{i})}{(1-tq^{n-i})}
\notag\\
&+\sum_{j=1}^{[(n-1)/2]}X^{n-2j}\,
\frac{(1-tq^{j})}{(1-tq^{n-j})}\prod_{i=0}^{j-1}
\frac{(1-q^{n-i})}{(1-q^{1+i})}\,\frac{(1-tq^{i})}{(1-tq^{n-i})},\\
\label{exacte+x}
E_{n}=X^{n}+&
\sum_{j=1}^{[n/2]}X^{2j-n}\,q^{n-j}\,
\frac{(1 - q^{j})}{(1\!-\!q^{n-j})}\prod_{i=0}^{j-1}
\frac{(1\!-\!q^{n-i-1})}{(1-q^{1+i})}
\frac{(1-tq^{i}\ ) }{(1\!-\!tq^{n-i-1})}
\notag\\
&+\sum_{j=1}^{[(n-1)/2]}X^{n-2j}\,q^j\,
\prod_{i=0}^{j-1}\frac{(1-q^{n-i-1})}{(1-q^{1+i})}\,
\frac{(1-tq^{i}\ ) }{(1-tq^{n-i-1})}.
\end{align}

These are  formulas (3.10 - 3.12) from \cite{C102} and \cite{ChO1}.
One can present them as infinite series, which are actually
finite and will terminate automatically. This is important for
what will follow. We note that T.~Koornwinder was the first to
explicitly calculate the formulas for the $E$\~polynomials for $A_1$.
One has:

\begin{align}\label{exacte-}
&E_{-n}\,=\,\sum_{j=0}^{\infty}X^{n-2j}\,
\frac{(1-tq^{j})}{(1-tq^{n-j})}\prod_{i=0}^{j-1}
\frac{(1-q^{n-i})}{(1-q^{1+i})}\,\frac{(1-tq^{i})}{(1-tq^{n-i})},\\
\label{exacte+}
&E_{n}=\sum_{j=0}^{\infty}X^{n-2j}\,q^j\,
\prod_{i=0}^{j-1}\,\frac{(1\!-\!q^{n-i-1})}{(1-q^{1+i})}\,
\frac{(1-tq^{i}\ ) }{(1\!-\!tq^{n-i-1})}\ \,  (n\!>\!0),
\end{align}
where the actual summation will be till $j\!=\!n\,$ in the
first formula (which works for $n=0$)
and till $j\!=\!n\!-\!1\,$ in the second;
the products there will vanish otherwise.
\smallskip

Setting\, $\th(X)\equal
\sum_{n=-\infty}^\infty q^{n^2/4}X^n$, the {\em global nonsymmetric
function} for $A_1$ is defined in \cite{C5}
as follows:
\begin{align}\label{gxla}
\frac{\th(X)\th(\La)}{\th(t^{1/2})}
G(X;\La)\, =\ &\sum_{n=-\infty}^\infty\,q^{\frac{\,\,n^2}{4}}\,
t^{\frac{|n|}{2}}
\,\frac{E_n^*(X)E_n(\La)}{\lr E_n E_n^*\mu_\circ\rr} \\
=\Psi(X,\La)\equal 
&\sum_{n=-\infty}^\infty\,q^{\frac{\,\,n^2}{4}}\,t^{\frac{|n|}{2}}
\,\frac{E_n(X)E_n^*(\La)}{\lr E_n E_n^*\mu_\circ\rr}\,. \notag
\end{align}

The following properties of the $G$\~function
are from  \cite{C5} (for any reduced root systems): 
\begin{align}\label{relyxt}
&Y(G)=\La^{-1}\,G,\
X^{-1}\,G=Y_\La(G),\
T(G)=T_\La(G),\\
\label{relyxta}
&G(X,\La=q^{n_\#})\ =\ \prod_{i=1}^\infty\,\frac{1-tq^i}{1-q^i}\, 
\frac{E_n(X)}
{E_n(t^{-1/2})}\,,\ \, n\in \Z,\\
&E_{n}(t^{-\frac{1}{2}})\!=\!
t^{-\frac{|n|}{2}}\!\!\!\prod_{0<j<|\tilde{n}|}
\frac{1\!-\!q^{j}t^{2}}{1\!-\!q^{j}t},\ 
|\tilde{n}|\!=\!\left\{\!\!\begin{array}{ccc}|n|\!+\!1
&\hbox{\!\!for\!\!}  & \!n\le 0, \\n &\hbox{\!\!for\!\!}  & 
\!n\!>\!0\end{array}\!\!\right\}.
\end{align}
Here $Y_\La$ and $T_\La$ are $Y,T$ where $X$ is replaced by $\La$. 
The relations from  (\ref{relyxt}), 
the Shintani-type formula (\ref{relyxta}) 
and the $X\leftrightarrow \La$ symmetry in  (\ref{gxla}) are the
key ingredients of our approach, extending that in \cite{ChW}.

The symmetric Macdonald polynomials $P_n (n\ge 0)$ for $A_1$, 
which are the Rogers polynomials, are given by the
formulas $P_n=E_{-n}+\frac{t-tq^n}{1-tq^n}E_n$. Accordingly,
\begin{align}\label{P=E+E}
(1+t)&\,\frac{P_n(X)}{P_n(t^{1/2})}\!=\!
\frac{1-t^2q^n}{1-tq^n}\frac{E_{-n}(X)}{E_{-n}(t^{-1/2})}
\!+\!\frac{t-tq^n}{1-tq^n}\frac{E_{n}(X)}{E_{n}(t^{-1/2})}
\hbox{\, for\, } n\!\ge\!0,
\notag\\
&(1+t)F(X,\La)\,=\,
\frac{t-\La^{-2}}{1-\La^{-2}}\,G(X,\La^{-1})
+\frac{t-\La^2}{1-\La^2}\,G(X,\La) \hbox{\, \,for }\\
&\frac{\th(\!X\!)\th(\!\La\!)}{\th(t^{1/2})}\,F(X,\La)\ =\ 
\Phi(X,\La)\,\equal\,
\sum_{n=0}^\infty\,q^{\frac{\,\,n^2}{4}}\,
t^{\frac{n}{2}}
\,\frac{P_n(X)P_n(\La)}{\lr P_n P_n\mu_\circ\rr}\notag
\end{align} 
apart from the zeros of $\th(\!X\!)\th(\!\La\!)$.
See the last three formulas of \cite{C5}.

\begin{theorem}\label{HCHNEXP} For the
function $G(X,\La)$ from (\ref{gxla}), let us assume
that $|q|<1$  and $|X|>|t|^{-1/2}|q|^{1/2}$. For
$\lan \mu\ran$ from (\ref{mucontone}),
\begin{align}\label{Sphexpan}
\Psi(X,\La)\ =&\,\ \lan \mu\ran\, \Bigl(
\si(\hbox{\small $\frac{1}{\La}$})\,
\th(X\!\La t^{\frac{1}{2}})\,\Xi_-(X,\La)\,
+\,\si(\La)
\,\th(\hbox{\small $\frac{X}{\La}$}t^{\frac{1}{2}})\,
\Xi_+(X,\hbox{\small $\frac{1}{\La}$})\Bigr),\notag
\\
\Xi_-(X,\La)\!=& 
\frac{1-t}{1\!\!-\!t\La^{\!-2}}\!+\!\!\sum_{j=1}^\infty
\!\frac{1-t q^j}{1\!-\!t\La^{\!-2}}
\left(\frac{q}{t}\right)^{j}\!X^{-2j}
\prod_{s=1}^j\frac{(1\!-\!tq^{s\!-\!1})
(1\!\!-\!tq^{s\!-\!1}\La^{\!2})}{(1-q^s)(1-q^s\La^{\!2})},
\notag\\
\Xi_+(X,\La)\!=&\, 
1+\sum_{j=1}^\infty
\frac{1-tq^{j}\,\La^{2}}{1-t\,\La^{2}}
\left(\frac{q}{t}\right)^{j}\!X^{-2j}\,
\prod_{s=1}^j\frac{(1\!-\!tq^{s-1})
(1\!-\!tq^{s-1}\,\La^{2})}{(1\!-\!q^s)(1-q^s\,\La^{2})},
\notag\\
&\hbox{where\, \,}
\si(\La)\,=\,\prod_{j=0}^\infty
\frac{1-tq^j \La^2}{1-q^j \La^2}=1+
\sum_{j=1}\La^{2j}\prod_{s=1}^j\frac{1-tq^{s-1}}
{1-q^{s}},
\end{align}
which is a $\,q,t$\~generalization of the Harish-Chandra 
$\,c$\~function. Using (\ref{P=E+E}), we arrive at
Theorem 2.3 from \cite{ChO1} for $X\mapsto X^{-1}:$
\begin{align}
&\Phi(X,\La)\!=\!
\lan \mu\ran\si(\hbox{\small $\frac{1}{\La}$})\,
\frac{\th(X\La t^{\frac{1}{2}})}{1+t}
\!\sum_{j=0}^\infty\!
\left(\frac{q}{t}\right)^{j}\!X^{-2j}\prod_{s=1}^j
\frac{(1\!-\!tq^{s\!-\!1})
(1\!-\!tq^{s\!-\!1}\La^{2})}{(1-q^s)(1-q^s\La^{2})}
\notag\\
+&\,\lan \mu\ran\,\si(\La)\,\frac{\th(X\La^{-1}t^{1/2})}{1+t}\,
\sum_{j=0}^\infty
\left(\frac{q}{t}\right)^{j}X^{-2j}\,\prod_{s=1}^j\frac{(1-tq^{s-1})
(1-tq^{s-1}\La^{-2})}{(1-q^s)(1-q^s\La^{-2})}.
\notag
\end{align}
\end{theorem}
{\it Proof.}
Let us present $X^{-n}E_{\pm n}$ 
in terms of $\La$. Setting
\begin{align}\label{e-exp}
&\tilde{\Xi}_+(X,\La)\!\equal\Xi_+(X,\La),\ \,\,
\tilde{\Xi}_-(X,\La)\!\equal\Xi_-(X,\La)\,
\frac{1\!-\!t\La^{\!-2}}{1\!-\La^{\!-2}}\\
=&\,\frac{1-t}{1\!-\!\La^{\!-2}}\!+\!\!\sum_{j=1}^\infty
\!\frac{1-t q^j}{1\!-\!\La^{\!-2}}
\left(\frac{q}{t}\right)^{j}\!X^{-2j}
\prod_{s=1}^j\frac{(1\!-\!tq^{s\!-\!1})
(1\!-\!tq^{s\!-\!1}\La^{\!2})}{(1-q^s)(1-q^s\La^{\!2})},\notag
\end{align}
we claim that 
\begin{align}\label{tildeXi}
&\frac{E_{n+1}}{X^{n+1}}\!=\!
\tilde{\Xi}_+(X,\La^{\!2}\!=\!\frac{1}{tq^{n+1}}),\ \,
\frac{E_{-n}}{X^n}\!=\!
\tilde{\Xi}_-(X,\La^{\!2}\!=\!\frac{1}{tq^{n}})
\hbox{\, for \,} n\!\ge\!0.
\end{align}
These two formulas are direct from
(\ref{exacte-},\ref{exacte+}), as well as
\begin{align}\label{e-expp}
&X^nE_{-n}\, =\, \hat{\Xi}_+(X,\La^{\!2}\!=\!\frac{1}{tq^{n}}) 
\hbox{\,\, for\, } n \ge 0, \hbox{\,\, where} \notag\\
&\hat{\Xi}_+(X,\La)\equal \sum_{j=0}^\infty
\!\frac{X^{2j}}{t^j}
\prod_{s=1}^j\frac{(1\!-\!tq^{s\!-\!1})
(1\!-\!tq^{s\!-\!1}\La^{\!2})}{(1-q^s)(1-q^s\La^{\!2})}.
\end{align}

Following \cite{C5} and using $Y(G)=\La^{-1}\,G$
from (\ref{relyxt}) (instead of the Macdonald eigenvalue
problem there), we obtain that 
\begin{align}\label{albeexp}
\frac{\Psi(X,\!\La)}{\lan\mu\ran}\!=\!
\al(\frac{1}{\La})\,\th(X\!\La t^{\frac{1}{2}})\,
\Xi_-(X,\!\La)\!+\!
\be(\La)\,\th(\frac{X}{\La} 
t^{\frac{1}{2}})\,\Xi_+(X,\!\frac{1}{\La})
\end{align}
for certain functions $\al(\La),\be(\La)$, which
do not depend on $X$ (see the general case below).
They can be determined as follows.

Let us find $\al(\La)$ and $\be(\La)$. We will use
the relation $\th(Xq^{m/2})=X^{-m}q^{-m^2/4}\,\th(X)$
and  the $\,\frac{1}{2}\Z$\~periodicity of 
$\,q^{x^2}\,\th(X)$ in terms of $x$ defined from $X=q^x$;
the latter periodicity readily implies the former relation.
We set $\La_n\equal
\!t^{\frac{1}{2}}q^{\frac{n}{2}}$ for $n>0$ and
$\La_{-n}\equal
(\!t^{\frac{1}{2}}q^{\frac{n}{2}})^{-1}$ for $n\ge 0$.

Applying the Shintani-type formula (\ref{relyxta}) for 
$E_{-n} (n\in \Z_+)$, let us establish the theorem when 
$\La=\La_{-n}$. One has: 
\begin{align}
\prod_{i=1}^{n}&
\frac{1\!-\!t^{2}q^{i}}{1\!-\!q^{i}t}\, 
\prod_{i=1}^\infty\,\frac{1-q^i}{1-tq^i}
\,t^{-\frac{n}{2}}\,G(X,\La_{-n}\!=\!(tq^n)^{-\frac{1}{2}})= 
E_{-n}(X)\, \Rightarrow\\
\!\prod_{i=1}^{n}
\frac{1\!-\!t^{2}q^{i}}{1\!-\!q^{i}t}\, 
&\prod_{i=1}^\infty\,\frac{1\!-\!q^i}{1\!-\!tq^i}\, 
\frac{\th(X)\th(t^{-\frac{1}{2}}q^{-\frac{n}{2}})}
{\th(t^{1/2})}G\bigl(X,(tq^n)^{-\frac{1}{2}}\bigr)\!=\! 
\th(X\!\La_{-n} t^{\frac{1}{2}})\frac{E_{-n}}{X^n},\notag\\
\!\prod_{i=1}^{n}
\frac{1\!-\!t^{2}q^{i}}{1\!-\!q^{i}t}\, 
&\prod_{i=1}^\infty\,\frac{1\!-\!q^i}{1\!-\!tq^i}\, 
\Psi\bigl(X, (tq^n)^{-\frac{1}{2}}\bigr)\ =\  
\th(X\!\La_{-n} t^{\frac{1}{2}})\,\frac{E_{-n}(X)}{X^n}.\notag
\end{align}
The last equality is exactly the decomposition for 
$\Psi(X,\La\!=\!\La_n)$ from (\ref{Sphexpan}). Indeed, 
using the product formula for $\lan \mu\ran$ from (\ref{mucontone})
and upon the division of the last equality by 
$\si(\La_n)\th(X\!\La_{-n} t^{\frac{1}{2}})$,
\begin{align*}
&\frac{\tilde{\Xi}_-(X,\La_{-n})}{\si(\La_{n})}=\lan\mu\ran 
\prod_{i=1}^{n}
\frac{1\!-\!t^{2}q^{i}}{1\!-\!q^{i}t}\, 
\prod_{i=1}^\infty\,\frac{1\!-\!q^i}{1\!-\!tq^i}\,
\Xi_-(X,\La_{-n})\\
=&\ \,\prod_{i=1}^\infty \frac{(1-tq^i)^2}{(1-t^2q^i)(1-q^i)}
\prod_{i=1}^{n}
\frac{1\!-\!t^{2}q^{i}}{1\!-\!tq^{i}}\, 
\prod_{i=1}^\infty\,\frac{1\!-\!q^i}{1\!-\!tq^i}\,\Xi_-(X,\La_{-n})\\
=&\prod_{i=n+1}^{\infty}
\frac{1- tq^i}{1\!-\!t^{2}q^{i}}\,\,\Xi_-(X,\La_{-n})
\,=\,\frac{1\!-\!t^{2}q^{n}}{1- tq^n}\,\,
\frac{\Xi_-(X,\La_{-n})}{\si(\La_n)},
\end{align*}
which is the relation between $\,\Xi\,$ and 
$\,\tilde{\Xi}\,$
from (\ref{e-exp}).


Thus $\al(\La_n)=
\si(\La_n)$ and $\be(\La_{-n})=0$ for $n\in \Z_+$.
Similarly, one checks that 
$\al(\La_{-n})=0$ and $\be(\La_{n})=
\si(\La_{n})$ for $n>0$. This is actually sufficient to determine
the coefficients $\al,\be$ uniquely by analyticity
considerations, which can be also established as follows.
\smallskip

{\sf Classical method.}
Let assume that $|X|$ approaches $\infty$ in (\ref{albeexp})
provided $|\La|>1$ and that the term with $\be$
does not contribute to the resulting {\em asymptotic\,} 
expansion. This is the classical track; 
we set $X=q^x, \La=q^{\la}$ providing  $\Re(x),\Re(\la)<0$ 
and $\Re(x\la)>0$. Under these assumptions, let
$\,\Psi^{-}_{-}(X,\La)\equal \lan \mu\ran^{-1}
\Psi(X,\!\La)\,\th(X\!\La t^{\frac{1}{2}})^{-1}$.
Then $\Psi^{-}_{-}(X,\La)$ 
approaches $\al(\La^{-1})(1-t)$ as $|X|\to\infty.$  
On the other hand, this limit  can be  
readily calculated using the symmetry $X\leftrightarrow \La$.
Namely, 
\begin{align}
\lim_{|X|\to \infty}& \Psi^{-}_{-}(X,\La)\!=\!
\lim_{n\to \infty} \Psi^{-}_{-}(\La, X_{-n}\!=\!t^{-1/2}q^{-n/2})=
\tilde{\Xi}_-(\La,|X|\!\to\!\infty)\notag\\
\label{findalbe}
&=\,1\!-\!t+\!\sum_{j=1}^\infty
(1\!-\!t q^j)\La^{-2j}
\prod_{s=1}^j\frac{1\!-\!tq^{s\!-\!1}}{1-q^s}
=(1\!-\!t)\si(\La^{-1}).
\end{align} 
Use the formula for $\tilde{\Xi}(X,\La)$ 
and the expansion of $\si$
from (\ref{Sphexpan}). 
\smallskip

Similarly, sending $|X|\to \infty\,$ in 
(\ref{albeexp}) for $\,|\La|<1\,$  provides the formula for the
$\be$\~coefficient. Now $\Re(\la)<0>\Re(x)$
and we the relation  $\Re (x\la)<0$ is imposed.
We set
\begin{align}\label{psipm}
\Psi^{+}_{-}(X,\La)\equal
\frac{\Psi(X,\!\La)}
{\lan\mu\ran}\, \th\!\hbox{\small 
$\left(\!\frac{\La}{X} t^{\frac{1}{2}}\!\right)^{\!-1}$}\!,\ 
\Psi^{-}_{+}\equal
\frac{\Psi(X,\!\La)}
{\lan\mu\ran}\, \th\!\hbox{\small 
$\left(\!\frac{X}{\La} t^{\frac{1}{2}}\!\right)^{\!-1}$}\!,
\end{align}
which are proportional to the expansions of $X^{n}E_{-n}$ and
$X^{-n}E_{n}$ as $\La=(tq^n)^{\pm 1/2}$ in terms of
$X^m$ for $m\in \pm \Z_+$ (with explicit coefficients
of proportionality due to the Shintani-type formulas). Parallel to 
the above consideration of $\Psi^{-}_{-}$\,, 
\begin{align}\label{belim}
\lim_{|X|\to \infty}
\Psi^{-}_{+}(X,\! \La)\!=\! 
\be(\La)\,\Xi_+(X\!\to\!\infty,\frac{1}{\La}) =
\be(\La).
\end{align}

Then we employ the $X\leftrightarrow \La$
symmetry and use (\ref{e-expp}):
\begin{align*}
\be(\La)\,&=
\lim_{|X|\to \infty}
\Psi^{+}_{-}(\La, X)
\!=\! \lim_{n\to \infty}
\Psi^{+}_{-}(\La, X\!=\!X_{-n})\!=\! 
\hat{\Xi}_-(\La,X\!\to\!\infty)\\
&=\lim_{|X|\to\infty}\,
\sum_{j=0}^\infty \frac{\La^{2j}}{t^j}
\prod_{s=1}^j\frac{(1\!-\!tq^{s\!-\!1})
(1\!-\!tq^{s\!-\!1}X^{\!2})}{(1-q^s)(1-q^sX^{\!2})}
=\si(\La).
\end{align*}
\vskip -1cm 
\sq
\smallskip

We note that $\lan \mu\ran$ appears in this theorem
as a result of direct calculation based on the 
constant term and the evaluation 
identities for the nonsymmetric Macdonald polynomials and 
the Shintani-type formula for global functions. 
The presence of $\lan \mu\ran$ can be seen 
in a more conceptual way, but we will not discuss this
in the present paper. Note that the $X\leftrightarrow \La$\~symmetry 
of $G(X,\La)$, which is the key, has no counterpart 
in the differential Harish-Chandra theory (though holds in
the rational limit, the theory of nonsymmetric Hankel transform).
\smallskip

The approach we use can be generalized to arbitrary
root systems. Following \cite{ChW}, one can calculate
the expansions of $X^{u(b)}E_{w(b)}$ for the nonsymmetric
Macdonald polynomials for any $b\in P_-$ and $u,w\in W$
(see \cite{C101} and below). They will be in 
terms of $X_\al$ for the roots $\al\in u(R_+)$ and
will coincide with $\,\si(X)\,$ up to certain explicit 
factors. The theorem above was stated only for $|X|>1$.
Following (\ref{psipm}), the functions $\Xi_-,\Xi_+$
can be naturally denoted by $\Xi^{-}_{-}$ and 
$\Xi^{-}_{+}$. Then the expansions of $\Psi$ for $|X|<1$
will be in terms of  $\Xi^{+}_{\pm}$. The calculation
of these expansions is very similar, but we will omit 
this in the present paper. Let us go to the general theory.

\setcounter{equation}{0}
\section{\sc General theory}
The extension of Theorem \ref{HCHNEXP} to arbitrary
reduced root systems (and its proof) is relatively straightforward.
We follow \cite{C4} (concerning $E$\~polynomials and
their main properties), \cite{C5} (for the global functions), 
\cite{ChW} (the $\si$\~function) and \cite{C101}. See 
\cite{Sto} for the global functions in the $C^\vee C_n$\~case.

Let $R=\{\al\}   \subset \R^n$ be a root system of type
$A,\!B,\ldots,\!G_2$
with respect to a euclidean form $(z,z')$ on $\R^n
\ni z,z'$, $W$ the Weyl group 
generated by the reflections $s_\al$,
$R_{+}$ the set of positive  roots
corresponding to fixed simple 
roots $\al_1,...,\al_n;$ $R_-=-R_+$. 
The form is normalized
by the condition  $(\al,\al)=2$ for 
{\em short\,} roots. 
The root lattice and the weight lattice are:
\begin{align}
& Q=\oplus^n_{i=1}\Z \al_i \subset P=\oplus^n_{i=1}\Z \om_i,
\notag
\end{align}
where $\{\om_i\}$ are fundamental weights:
$ (\om_i,\al_j^\vee)=\de_{ij}$ for the
coroots $\al^\vee=2\al/(\al,\al).$
Replacing $\Z$ by $\Z_{+}=\{m\in\Z, m\ge 0\}$, we obtain
$Q_+, P_+;$ we will constantly use $P_-=-P_+$. Let
$\iota(b)=b^\iota=-w_0(b)$ for the element $w_0$ of
the maximal length in $W$.
\smallskip

Setting 
$\nu_\al\equal (\al,\al)/2$,
the vectors $\, \tal=[\al,\,\nu_\al j\,] \in
\R^n\times \R \subset \R^{n+1}$
for $\al \in R, j \in \Z $ form the
{\em twisted affine root system\,}
$\tR \supset R$ ($z\in \R^n$ are identified with $ [z,0]$).
The corresponding set
$\tR_+$ of positive roots is 
$R_+\cup \{[\al,\,\nu_\al j\,],\ \al\in R, \ j > 0\}$.
\smallskip

The {\em extended Weyl group\,} $ \hW$ is $W\lsmash P$, where
the corresponding action in $\R^{n+1}$ is 
\begin{align}
&(wb)([z,\ze])\ =\ [w(z),\ze-(z,b)] \for w\in W,\, b\in P.
\notag
\end{align}
The {\em length } in $\hW$ is defined as 
follows:
\begin{align*}
&l(\hw)=|\la(\hw)| \for \la(\hw)\equal\tR_+\cap \hw^{-1}(-\tR_+).
\end{align*}
\smallskip

For $\tal=[\al,\,\nu_\al j\,] \in \tR,\ 0\le i\le n$, we set
\begin{align*}
&   t_{\tal} =t_{\al}=t_{\nu_\al}=q_\al^{k_\nu} ,\ \, 
q_{\tal}=q^{\nu_\al}, \ \, t_i = t_{\al_i},\,  q_i=q_{\al_i},\\
 \rho_k\equal \frac{1}{2}\!&\sum_{\al>0} k_\al \al=
k_{\sht}\rho_{\sht}\!+\!k_{\lng}\rho_{\lng},\ \,
\rho_\nu=\frac{1}{2}\!\sum_{\nu_\al=\nu} \al=
\!\!\sum_{\nu_i=\nu,i>0}  \om_i,
\end{align*}
where \,{\em\small sht,\ lng\,} are used instead 
of $\nu$. 
\smallskip

Let $\v\,$ by the space of Laurent polynomials
in terms of $X_b (b\in P)$ satisfying the
multiplicative property $X_{b+c}=X_bX_c$;
the coefficients are taken from 
$\Z[q^{\pm 1/m},t_\nu^{\pm 1/2}]$,
where $(P,P)\in \frac{1}{m}\Z$. 
For any $b\in P$ and generic $q,t$,
the {\em nonsymmetric Macdonald polynomials\,}
are uniquely determined by the relations
\begin{align}\label{yaebdef}
&Y_a(E_b)\,=\,q^{-(a,b-u_b(\rho_k))}E_b \hbox{\,\, and 
the coefficient of\,\,} 
X_b\in E_b \hbox{\, is \,} 1\\
&\hbox{for\, } u_b\in W \hbox{\, of 
minimal possible length such that\, } u_b(b)\in P_-,\notag
\end{align}
where $Y_a (a\in P)$ are the
{\em difference Dunkl operators\,}; they act in $\v$
and are pairwise commutative. See Theorem 4.1 from \cite{C4} and 
Proposition 3.3.1 from \cite{C101}. The coefficients of $E_b$ 
belong to $\Q(q,t_\nu)$.
Setting $b_\#\equal b-u_b(\rho_k)$, we can rewrite  (\ref{yaebdef})
as 
\begin{align}\label{yebx}
Y_a(E_b)\,=\,X_b^{-1}(q^{b_\#})E_b,\hbox{\, where\, }
X_a(q^b)\equal q^{(a,b)} \for a,b\in P. 
\end{align}

Let $\xi(\al)\equal 0,1$ respectively for $\al\in R_+,R_-$. 
The Main Theorem in \cite{C4} 
and formula (3.3.16) in \cite{C101} state that
\begin{align}\label{macde-eval}
&E_{b}(q^{-\rho_k})=
q^{(\rho_k,b_-)}
\prod_{\al>0}\prod_{j=1}^{j(b,\al)}
\Bigl(
\frac{
1- q_\al^{j}t_\al X_\al(q^{\rho_k})
 }{
1- q_\al^{j}X_\al(q^{\rho_k})
}
\Bigr) \for b\in P,\\
&\hbox{where\, }
b_-\equal u_b(b)\!\in\! P_-,\ \,j(b,\al)=-(\al^{\!\vee},b_-)-
\xi(u_b^{-1}(\al)).\notag
\end{align}
\smallskip

We will also need 
\begin{align}\label{mudefin}
&\mu(X;q,t)\!\equal\!\!\prod_{\al \in R_+}
\prod_{j=0}^\infty \frac{(1\!-\!X_\al q_\al^{j})
(1\!-\!X_\al^{-1}q_\al^{j+1})
}{
(1\!-\!X_\al t_\al q_\al^{j})
(1\!-\!X_\al^{-1}t_\al^{}q_\al^{j+1})},\\
&\langle\mu\rangle\ =\ \prod_{\al \in R_+}
\prod_{i=1}^{\infty} \frac{ (1- q^{(\rho_k,\al)+i\nu_\al})^2
}{
(1-t_\al q^{(\rho_k,\al)+i\nu_\al})
(1-t_\al^{-1}q^{(\rho_k,\al)+i\nu_\al})
},
\label{consterm}
\end{align}
where $\lan \,\cdot\,\ran$ is the constant term functional
(the coefficient of $X^0$).
See e.g. (3.3.1),(3.3.2) in \cite{C101}. We will use below
the norms
$\lan E_b E_b^\star\mu_\circ\ran$ for 
$\mu_\circ\equal\mu/\lan \mu\ran$ from (3.4.2) there, but do not 
need the exact formulas for them.
The conjugation $\star\,$ is 
as follows: 
$X_b^\star=X_b^{-1}, q^\star=q^{-1}, t_\nu^\star=t_\nu^{-1}$.
\smallskip

We will omit the definition of the symmetric Macdonald
polynomials $P_b (b\in P_-)$ in this paper; see
(3.3.12) and (3.3.14) from \cite{C101} and (3.4.3) there
for the norms $\lan P_b P_{b^\iota}\mu_\circ\ran$. 
Recall that $b^\iota=-w_0(b)$. Formula (3.3.23) from
\cite{C101} states:
\begin{align*}
P_{b_-}(q^{-\rho_k})=P_{b_-}(q^{\rho_k})=
E_{b_+}(q^{-\rho_k})
\prod_{\al>0}
\Bigl(
\frac{1- t_\al X_\al(q^{\rho_k})}{1- X_\al(q^{\rho_k})}
\Bigr). 
\end{align*}
The following formula
is from (4.13) \cite{C4} and (3.3.15) from \cite{C101}:
\begin{align}
&P_{b_-}\ =\
\sum_{c\in W(b_+)}
\prod_{(\al,c)>0}\frac{t_\al-X_\al(q^{c_\#})}
{1-X_\al(q^{c_\#})}E_c,\ \al\in R_+.
\label{symmviae}
\end{align}
See \cite{M4} for $k_{\sht}=k_{\lng}\in\Z_+$.
We need its variant in terms of 
\begin{align*}
&E_b'\equal E_{b}/E_b(q^{-\rho_k}) \hbox{\, and\, } 
P_{b_-}''\equal \p_{\!R}(t)
P_{b_-}/P_{b_-}(q^{-\rho_k}),\\
&\hbox{where\, }
\p_{\!R}(t)=\prod_{\al>0}\frac{1-t_\al q^{(\al,\rho_k)}}
{1-q^{(\al,\rho_k)}}
\hbox{\, is the Poincar\'e polynomial}.
\end{align*}
 For $\al\in R_+$,
\begin{align}\label{symmvias}
&P''_{b_-}\!=\!\!\!\sum_{c\in W(b_+)}
\,\prod_{u_c^{-1}(\al)>0}\!
\Bigl(\frac{1- t_\al X_\al(q^{\rho_k-b_-})}
{1\!-\!X_\al(q^{\rho_k-b_-})}\Bigr)\!
\prod_{(\al,c)>0}\!
\frac{t_\al\!-\!X_\al(q^{c_\#})}{1-X_\al(q^{c_\#})}E'_c.
\end{align}
Switching in the first product to $\be=u_c^{-1}(\al)$,
the condition $u_c^{-1}(\al)>0$ for $\al>0$
becomes $(\be,c)\le 0$ for $\be>0$. Also $\rho_k-b_-=
-u_c(c_\#)$.
Therefore for $\al,\be>0$,
\begin{align}\label{symmviass}
P''_{b_-}\!=\!\!\!&\sum_{c\in W(b_+)}
\,\prod_{(\be,c)\le 0}\!
\Bigl(\frac{1- t_\be X^{-1}_\be(q^{c_\#})}
{1\!-\!X_\be^{-1}(q^{c_\#})}\Bigr)\!
\prod_{(\al,c)>0}\!
\frac{t_\al\!-\!X_\al(q^{c_\#})}{1-X_\al(q^{c_\#})}E'_c\\
= &\sum_{c\in W(b_+)}
\,\prod_{\al>0}\!
\frac{t_\al\!-\!X_\al(q^{c_\#})}{1-X_\al(q^{c_\#})}E'_c.\notag
\end{align}

\medskip

The following construction is from Theorem 5.4 and 
Corollary 7.3 of \cite{C5}. The function $G(X,\La)$
introduced in (\ref{gexla}) below is called 
{\em global nonsymmetric $q,t$\~spherical function}.
We will use another set of (pairwise commutative)
variables $\La_b$. Sometimes it will be convenient
to put $\,X=q^x$. Then $\,X_b=q^{(x,b)}$ for $b\in P$. 
Accordingly, $w(X_b)=X_{w(b)}$ for $w\in W$ and 
$w\bigl(f(X)\bigr)=f(w^{-1}(X))=f(q^{w^{-1}(x)})$ for any function
$f$ of $X$ (notice $w^{-1}$). The same notations will be
used for $\La=q^{\la}$; we put $w_\La$ instead of $w$ if
this action can be confused with that in terms of $X$.

We set $\th(X)\equal\sum_{b\in P} q^{-(b,b)/2}X_b$;
obviously $w(\th(X))=\th(X)=th(X^{-1})$ for $w\in W$
and $q^{\frac{(x,x)}{2}}\th(X)$ is $P$\~periodic in
terms of $x$. If $|q|<1$, then $\th(X)$ is convergent
and holomorphic anywhere.
 
\begin{theorem}\label{GLOBNSSPH}
(i) The Laurent series
\begin{align}
\Psi(X,\La)=\Psi(X,\La;q,t)&\equal
\sum_{b\in B} q^{(b,b)/2-(b_-,\rho_k)} \
\frac{E_b^\star(X) \, E_b(\La)}
{\lan E_b E_b^\star\mu_\circ \ran}
\label{psiexla}
\end{align}
is well defined with coefficients in
$\Q[t][[q^{\frac{1}{2m}}]]$.
For $|q|<1$, $\Psi$ converges to an entire function
of $X,\La$, provided $t_\nu$ are chosen so that all $E$\~polynomials
exist (the conditions $|t_\nu|<1$ are sufficient).
Accordingly, 
\begin{align}
G(X,\La) \equal 
\frac{\th(q^{\rho_k})}{\th(X)\th(\La)}\Psi(X,\La;q,t)
\label{gexla}
\end{align}
is a meromorphic function of
$X,\La$, which is analytic
apart from the zeros of $\th(X)\th(\La)$.


(ii) The function $G(X,\La)$ satisfies the relations
\begin{align}\label{GsymT}
G(X,\La)=G(\La,X),\ \,
Y_a(G(X,\La))=\La_{a}^{-1}G(X,\La) \for a\in P.
\end{align}
For an arbitrary $b\in P$, one has the Shintani-type
formulas:
\begin{align}\label{Shinqt}
& G(X, q^{b_\#})\ =\ 
\frac{E_b(X)}{E_b(q^{-\rho_k})} 
\prod_{\al\in R_+}\prod_{ j=1}^{\infty}\Bigl(\frac{ 
1- q^{(\rho_k,\al)+\nu_\al j}}{
1-t_\al^{-1}q^{(\rho_k,\al)+\nu_\al j} }\Bigr). 
\end{align}

(iii) Let us define 
the symmetric global function
 $F(X,\La)$  via
\begin{align}\label{fpxla}
\frac{\th(X)\th(\La)}{\th(q^{\rho_k})}
F(X,\!\La)\!=\!
\Phi(X,\!\La)\!\equal\!
\sum_{b\in P_-} q^{\frac{(b,b)}{2}-
\!\frac{(b,\rho_k)}{2}} \
\frac{P_b(X) P_{b^\iota}(\La)}
{\lan P_b P_{b^\iota}\mu_\circ \ran}.
\end{align}
Then (\ref{Shinqt}) holds 
with the same coefficient
of proportionality:
\begin{align}\label{Shinqtsym}
F(X, q^{b -\rho_k})\!&=\! 
\frac{P_b(X)}{P_b(q^{\rho_k})} 
\prod_{\al\in R_+}\prod_{ j=1}^{\infty}\Bigl(\frac{ 
1- q^{(\rho_k,\al)+\nu_\al j}}{
1\!-\!t_\al^{-1}q^{(\rho_k,\al)+\nu_\al j} }\Bigr)
\hbox{\, for\, } b=b_-,\\
\hbox{and (\ref{symmviass}) implies\, \,}
&\p_{\!R}(t)\,F(X,\La)\,=\,\sum_{w\in W}
\,\,\prod_{\al>0}
\frac{t_\al\!-\!\La_{w(\al)}}{1-\La_{w(\al)}}\,\,G(X,w(\La)).
\notag
\end{align}
\sq
\end{theorem}

The proof is based on the fundamental
fact that $G(X,\La)$ represents the
Fourier transform of DAHA. For instance,
the Shintani-type formula (\ref{Shinqt}) follows
from formula (3.4.13) in Theorem 3.4.2 \cite{C101},
which states that the Fourier-images of $E_b/E_b(q^{-\rho_k})$ 
(with respect to the DAHA-automorphism
$\vep$) are the corresponding delta-functions.

\smallskip

\begin{theorem}\label{HCHGEN} 
(i) {\sf  $\mathbf {\La}$\~Stabilization.\,}
For any given $w\in W$, there exists a unique series 
$\Xi^{(w)}(X,\La)=\sum_{a\in Q_+} A^w_a(\La)X_{a}^{-1}$
with the coefficients that are rational functions in terms
of $q,t_\nu,\La_\al$ such that
\begin{align}\label{asymser} 
X_{b^\iota_-}E_{b}\,=\, \Xi^{(w)}\bigl(
X,\La\!=\!q^{b_- -\rho_k}\bigr) \for w=u_b^{-1},\, b\in P;
\hbox{\, see (\ref{yaebdef})}.
\end{align} 
The functions $A^w_a(\La)$ are regular if
$\,\La_\al\not\in q_\al^{\Z}\,$ for every $\,\al\in R_+$
and have formal expansions in the ring 
$\,\Q[t_\nu^{\pm 1}]\,[[\,q_\al\La_\al,\,\La_\al^{-1}(\al>0)\,]].$ 
Assuming $|q|<1$, let $X=q^x$. There exists a constant\,
$C\!>\!0$ (dependent on $q,t_\nu$) such that the
series $\Xi^{(w)}$  converges for any $\,w\in W$
when $\,(\Re(x),\al)<\!-C\,$ for every 
$\al\in R_+$ and $\,\La$ can be arbitrary apart from 
the singularities of $A_a^w$.
\smallskip

(ii) {\sf Asymptotic decomposition.}
Let  $G(X,\La),\,\Psi(X,\La)$ be from (\ref{gexla})
and $\lan \mu\ran$ from (\ref{consterm}). Recall that
$\,w_\La(f(\La))=f\bigl(w^{-1}(\La)\bigr)\,$ and 
$\,w_\La(\La_b)=\La_{w(b)}\,$ for $\,w\in W, b\in P$; also,
$\La^\iota=-w_0(\La),\ (\La^\iota)_b=\La_{b^\iota}$.
\begin{align}\label{SphPsi}
&\hbox{Setting\,\, }
\si_\ast(\La)\equal\,\prod_{\al>0}\prod_{j=1}^\infty
\frac{1-t_\al q_\al^j \La^{-1}_\al}{1-q_\al^j \La^{-1}_\al},\ \, 
\Psi(X,\La)\,=\\
=\lan \mu\ran\sum_{w\in W}&
w_\La\Bigl(\si_\ast(\La)\,
\th\bigl(\La^\iota X q^{\rho_k}\bigr)
\Xi^{(w)}(X,\La)\Bigr)\prod_{\al>0>w(\al)}
\Bigl(\frac{1\!-\!t_\al\La^{-1}_{w(\al)}}
{1 - \La^{-1}_{w(\al)}}\Bigr),
\notag
\end{align}
provided the convergence of the right-hand side;
$\si_\ast(\La)$ is a modification of 
the $q,t$-generalization $\,\si\,$
of the Harish-Chandra  $c$\~function from formula
(4.5) in \cite{ChW} with $j\ge 1$ instead of
$\,j\ge 0\,$ and $\,\La\mapsto \La^{-1}$. 
\end{theorem}

{\em An outline of the proof}.
The stabilization in $(i)$ can be deduced directly
from the eigenvalue problem (\ref{yebx}), which is the
definition of $E_b$ in this paper. However it is more
convenient to use the intertwining operators from 
\cite{C1} (due to Knop and Sahi for $A_n$). 
One can also used the formulas from \cite{HHL,RY,OS}
based on the intertwining operators; see (the most
general) Theorem 3.13 in \cite{OS}. The intertwining
operators are not creation operators any longer (as
for the $E$\~polynomials), but provide recurrence
relations that are sufficient to justify $(i)$.
We will omit the details in this paper.
\smallskip

The key in $(ii)$ is the verification that the asymptotic 
expansion of $\,\Psi(X,\La)$ for 
$\,\La=q^{b_\#}=q^{u_b^{-1}(b_- -\rho_k)}\,$
coincides with the term with $w=u_b^{-1}$ in the right-hand side 
of (\ref{SphPsi}). Here $b=u_b^{-1}(b_-)$; see (\ref{yebx}). 
\smallskip 

First of all, let us check that the multiplier 
\begin{align}\label{siwla}
v\bigl(\si_\ast(\La)\bigr)\prod_{\al>0>v(\al)}
\Bigl(\frac{1\!-\!t_\al\La^{-1}_{v(\al)}}
{1 - \La^{-1}_{v(\al)}}\Bigr) \for v\in W
\end{align} 
is nonzero only when $v=u_b^{-1}$ 
for such $\La$.
Indeed, if $v\neq u_b^{-1}$, there exists $\al>0$
and a simple root $\al_i$ such that $u_b v(\al)=-\al_i$.
This results in 
$\,\bigl(u_b^{-1}(\rho)\,,\,v(\al)\bigr)=-1$ and
 $\bigl(u_b^{-1}(b_-)\,,\,v(\al)\bigr)=
(b_-\,, -\al_i)\ge 0$. If 
$(b_-\,,\al_i)=0$, then $u_b^{-1}(\al_i)>0$ 
(we use the minimality of $u_b$) and
$v(\al)<0$, which is the inequality from the product part
in (\ref{siwla}).
Therefore the term $\bigl(1-q^j t_\al \La_{v(\al)}^{-1}\bigr)$ 
in the numerator of (\ref{siwla}) will vanish at 
$\La=q^{b_\#}$ for $j=-(b_-\,,\al_i)\ge 0$.
\smallskip

We will use the relation 
$\,\th(Xq^b)=X_b^{-1}q^{-(b,b)/2}\th(X)\,$
for $b\in P$, which follows from the $P$\~periodicity
of $q^{\frac{(x,x)}{2}}\th(X)$ and readily results in
\begin{align*}
\th(X\!q^{b^\iota -\rho_k} q^{\rho_k})\,=\,
q^{-(b,\rho_k)}X_{b^\iota}^{-1}\,
\frac{\th(X)\th(q^{b -\rho_k})}{\th(q^{\rho_k})}.
\end{align*}

Let us consecutively apply formulas (\ref{macde-eval}), 
(\ref{Shinqt}) and (\ref{consterm}) (the Shintani 
formula, the evaluation and the constant term formulas).
We will begin with $b=b_-\in P_-, w=$id; the case of
general $b\in P$ is parallel.
\smallskip

We evaluate $\lan \mu\ran \si_\ast(q^{b-\rho_k})
\th(X\!q^{b^\iota -\rho_k} q^{\rho_k})
\Xi^{(\hbox{\tiny id})}\bigr(X,q^{b -\rho_k}\bigl)$,
which is the first and the only nonzero term in the decomposition
from (\ref{SphPsi}) for $\La=q^{b_\#}=q^{b-\rho_k}$.
Using $(i)$ for $b=b_-$,
\begin{align*}
&
\th(X q^{b^\iota -\rho_k} q^{\rho_k})\,
\Xi^{(\hbox{\tiny id})}\bigr(X,q^{b -\rho_k}\bigl)
=
q^{-(b,\rho_k)}X_{b^\iota}^{-1}\bigl(X_{b^\iota}E_b\bigr)
\frac{\th(X)\th(q^{b -\rho_k})}{\th(q^{\rho_k})}\\
=\,&
q^{-(b,\rho_k)}\Bigl(E_b(q^{-\rho_k})E'_b\Bigr)
\frac{\th(X)\th(q^{b -\rho_k})}{\th(q^{\rho_k})}
=
q^{-(b,\rho_k)}E_b(q^{-\rho_k})\\
\times&\Bigl(
\prod_{\al>0}\prod_{ j=1}^{\infty}\bigl(\frac
{1-t_\al^{-1}q^{(\rho_k,\al)+\nu_\al j}}
{1- q^{(\rho_k,\al)+\nu_\al j}}
\bigr)G(X,q^{b-\rho_k})\Bigr)
\frac{\th(X)\th(q^{b -\rho_k})}{\th(q^{\rho_k})}\\
=\,&
q^{-(b,\rho_k)}E_b(q^{-\rho_k})
\prod_{\al>0}\prod_{ j=1}^{\infty}\bigl(\frac
{1-t_\al^{-1}q^{(\rho_k,\al)+\nu_\al j}}
{1- q^{(\rho_k,\al)+\nu_\al j}}
\bigr)\Psi(X,q^{b-\rho_k})\\
=\,&
\prod_{\al>0}\!\left(\prod_{j=1}^{-(b,\al^{\!\vee})}
\frac{1\!-\! q_\al^{j}t_\al X_\al(q^{\rho_k})}
{1- q_\al^{j}X_\al(q^{\rho_k})}
\prod_{j=1}^{\infty}
\frac{1\!-\!t_\al^{-1}q^{(\rho_k,\al)\!+\!\nu_\al j}}
{1- q^{(\rho_k,\al)+\nu_\al j}}\right)\!
\Psi(X,q^{b-\rho_k}).
\end{align*}

Thus, 
\begin{align}
&\lan \mu\ran \si_\ast(q^{b_--\rho_k})
\th(X\!q^{b^\iota_- -\rho_k} q^{\rho_k})
\Xi^{(\hbox{\tiny id})}\bigr(X,q^{b_- -\rho_k}\bigl)
/\Psi(X,q^{b_--\rho_k})\notag
\\
=&\prod_{\al>0}\prod_{j=1}^{\infty} 
\frac{ (1- q_{\al}^j q^{(\rho_k,\al)})^2}
{(1-t_\al q_{\al}^j q^{(\rho_k,\al)})
(1-t_\al^{-1}q_{\al}^j q^{(\rho_k,\al)})} \notag\\
\label{Xievalmiddle}
\times&
\prod_{\al>0}\prod_{j=1}^\infty
\frac{1-t_\al q_\al^j q^{-(b_-,\al)+(\rho_k,\al)}}
{1-q_\al^j q^{-(b_-,\al)+(\rho_k,\al)}}\\
\times&
\prod_{\al>0}\!\left(\prod_{j=1}^{-(b_-,\al^{\!\vee})}
\Bigl(\frac{1- t_\al q_\al^{j} q^{(\rho_k,\al)}}
{1- q_\al^{j}q^{(\rho_k,\al)}}\Bigr)
\prod_{ j=1}^{\infty}\Bigl(\frac
{1-t_\al^{-1}q^{(\rho_k,\al)+\nu_\al j}}
{1- q^{(\rho_k,\al)+\nu_\al j}}
\Bigr)\!\right)\!=\!1.\notag
\end{align}
This cancelation readily results from
\begin{align*}
&\prod_{\al>0}\left(\prod_{j=1}^\infty
\frac{1-t_\al q_\al^j q^{-(b_-,\al)+(\rho_k,\al)}}
{1-q_\al^j q^{-(b_-,\al)+(\rho_k,\al)}}
\prod_{j=1}^{-(b_-,\al^{\!\vee})}
\Bigl(\frac{1- t_\al q_\al^{j} q^{(\rho_k,\al)}}
{1- q_\al^{j}q^{(\rho_k,\al)}}\Bigr)\right)=\\
&\prod_{\al>0}\prod_{j=1}^{\infty}
\Bigl(\frac{1- q_\al^{j}t_\al q^{(\rho_k,\al)}}
{1- q_\al^{j}q^{(\rho_k,\al)}}\Bigr).
\end{align*}

\comment{
\begin{align}\label{macde-evalx}
&E_{b}(q^{-\rho_k})=
q^{(\rho_k,b_-)}
\prod_{\al>0}\prod_{j=1}^{j(b,\al)}
\Bigl(
\frac{
1- q_\al^{j}t_\al X_\al(q^{\rho_k})
 }{
1- q_\al^{j}X_\al(q^{\rho_k})
}
\Bigr) \for b\in P,\\
&\hbox{where\, }
b_-\equal u_b(b)\!\in\! P_-,\ \,j(b,\al)=-(\al^{\!\vee},b_-)-
\xi(u_b^{-1}(\al)).\notag
\end{align}
\begin{align}\label{Shinqtx}
& G(X, q^{b_\#})\ =\ 
\frac{E_b(X)}{E_b(q^{-\rho_k})} 
\prod_{\al\in R_+}\prod_{ j=1}^{\infty}\Bigl(\frac{ 
1- q^{(\rho_k,\al)+\nu_\al j}}{
1-t_\al^{-1}q^{(\rho_k,\al)+\nu_\al j} }\Bigr). 
\end{align}
\begin{align}\label{mudefinx}
&\mu(X;q,t)\!\equal\!\!\prod_{\al \in R_+}
\prod_{j=0}^\infty \frac{(1\!-\!X_\al q_\al^{j})
(1\!-\!X_\al^{-1}q_\al^{j+1})
}{
(1\!-\!X_\al t_\al q_\al^{j})
(1\!-\!X_\al^{-1}t_\al^{}q_\al^{j+1})},\\
&\langle\mu\rangle\ =\ \prod_{\al \in R_+}
\prod_{i=1}^{\infty} \frac{ (1- q^{(\rho_k,\al)+i\nu_\al})^2
}{
(1-t_\al q^{(\rho_k,\al)+i\nu_\al})
(1-t_\al^{-1}q^{(\rho_k,\al)+i\nu_\al})
},
\label{constermx}
\end{align}
}

When $b\in P$ is arbitrary, let $w=u_b^{-1}\,$ and
$\La=u_b^{-1}(q^{-b_--\rho_k})$. Then  
\begin{align*}
\prod_{\al>0>w(\al)}
\frac{1\!-\!t_\al\La^{-1}_{w(\al)}}
{1 - \La^{-1}_{w(\al)}}=
\prod_{\al>0>w(\al)}
\frac{1-t_\al q^{-(b_-,\al)+(\rho_k,\al)}}
{1-q^{-(b_-,\al)+(\rho_k,\al)}}
\end{align*}
must be added to (\ref{Xievalmiddle}).
Respectively, the upper limit $\,-(b_-,\al^{\!\vee})\,$
for $j\,$ (the last line there) must be diminished
to $\,j(b,\al)=-(b_-,\al^{\!\vee})-
\xi(u_b^{-1}(\al))$, where $\xi(\al)=0,1$
for $\al\in R_{\pm}$; see (\ref{macde-eval}).
These adjustments complement each other and the
result will be $1$, as well as for $b=b_-$.
\smallskip

Next, we use the decomposition
\begin{align}\label{Psidecomp}
&\Psi(X,\La)\,=\lan \mu\ran\sum_{w\in W}
w_\La\Bigl(\varpi_w(\La)\,
\th\!\left(\La^\iota X q^{\rho_k}\right)
\Xi^{(w)}\bigl(X,\La\bigr)\Bigr),
\end{align}
where the coefficients $\varpi_w$ do not depend on $X$.
This is a general fact, which follows from
the difference Dunkl eigenvalue problem in 
(\ref{GsymT}): $Y_a(G(X,\La))=\La_{a}^{-1}G(X,\La)$ for $a\in P$.
Here you need to consider the latter in the {\em $W$\~spinors\,},
treating the action of $w\in W$ as the corresponding
permutation of the {\em independent\,} spinor components. 
Note that this is parallel to what was done in \cite{Op} in the 
differential setting.
The (formal) origin of the $W$\~spinors are  
\cite{C13} (where they were used for the Cherednik-Matsuo theorem) 
and \cite{Op}; see \cite{C102,ChO2} for the exact definitions and 
a comprehensive discussion.

Here generally, $\varpi_w$ can be $P$\~periodic functions
in terms of $x$. However in our setting, all 
functions are Laurent series in terms of $X$ and 
therefore $\,\varpi_w\,$ must be double periodic $x$\~functions
for $P\oplus \frac{2\pi \imath}{\log(q)}P$. 
Since they have no singularities in $x$ for sufficiently large 
negative $\Re(x)$ due to the convergence
of $\Psi$, they do not depend of $X$.
\smallskip

We conclude that for the Kronecker delta $\de_{v,w}$,
\begin{align}\label{varpicoef}
&w\bigl(\varpi_v(\La)\bigr)
=\de_{v,w}\,w\bigl(\si_\ast(\La)\bigr)\prod_{\al>0>w(\al)}
\Bigl(\frac{1\!-\!t_\al\La^{-1}_{w(\al)}}
{1 - \La^{-1}_{w(\al)}}\Bigr), \hbox{\, where\,}\\
&\La\!=\!w(q^{b_--\rho_k})\!=\!q^{b_\#} \for b\in P,\ \,
w\!=\!u_b^{-1},\ \,b_-\!=\!u_b(b)\in P_-.\notag
\end{align}

Relations from (\ref{varpicoef}) are actually sufficient to
fix $\{\varpi_w\}$ uniquely via the analyticity considerations.
A usual (classical) justification of similar facts is
what we did for $A_1$, using the $X\leftrightarrow \La$\~\,symmetry 
and considering the asymptotic sectors
where exactly one of the terms in the decomposition 
(\ref{SphPsi}) survives. One needs to know in this approach the
limits of $X_{w(b_-)}E_{b_-}$ as $|X_{w(\al)}|\to \infty$ for
all $\al>0$; cf. (\ref{findalbe}). These limits are 
relatively straightforward to find. As in \cite{ChW}, we obtain
the recurrence relations for them (coinciding with those for the
corresponding $\varpi_w$). We will omit the details in this paper.
\sq
\smallskip

{\sf An orbit-sum formula.}
There is an interesting application of
the decomposition theorem to the $E$\~polynomials, which is
another (``global") relation connecting $E_b$ and $\Xi^{(w)}$.
Using the $X\leftrightarrow\La$\~\,symmetry of $G(X,\La)$,
we can substitute $X=q^{b_--\rho_k}$ in (\ref{SphPsi}).
One has:
\begin{align}\label{SphPsix}
&\th(\La)^{-1}\,q^{\frac{(b,b)}{2}}\,\Psi(q^{b_--\rho_k},\La)\\
=\prod_{\al>0}\,&\left(\,\prod_{j=1}^{-(b_-,\al^{\!\vee})}
\Bigl(\frac{1- t_\al q_\al^{j} q^{(\rho_k,\al)}}
{1- q_\al^{j}q^{(\rho_k,\al)}}\Bigr)
\,\prod_{ j=1}^{\infty}\Bigl(\frac
{1-t_\al^{-1}q_\al^j q^{(\rho_k,\al)}}
{1- q_\al^j q^{(\rho_k,\al)}}\Bigr)\right)^{\!\!-1}\!
E_b(\La)\notag\\
=\!\lan \mu\ran \sum_{w\in W}&\,
w_\La\Bigl(\!\si_\ast(\La)\,
\frac{\th\bigl(\La^\iota q^{b_--\rho_k} q^{\rho_k}\bigr)}
{\th(\La)q^{-\frac{(b,b)}{2}}}
\Xi^{(w)}\bigl(q^{b_--\rho_k},\La\bigr)\!\Bigr)\!\!\!
\prod_{\al>0>w(\al)}\!
\frac{1\!-\!t_\al \La^{-1}_{w(\al)}}
{1 - \La^{-1}_{w(\al)}}\notag\\
=\!\lan \mu\ran\sum_{w\in W}&\,
w_\La\Bigl(\si_\ast(\La)\,\La_{b^\iota_-}^{-1}\,
\Xi^{(w)}\bigl(q^{b_--\rho_k},\La\bigr)\Bigr)\!
\prod_{\al>0>w(\al)}
\frac{1\!-\!t_\al \La^{-1}_{w(\al)}}
{1 - \La^{-1}_{w(\al)}}\,.\notag
\end{align}

Replacing $\lan \mu\ran$ by the corresponding product,
we obtain the following identity: 
\begin{align}\label{SphPsixx}
&\prod_{\al>0}\,\prod_{j=1}^{\infty}
\Bigl(\frac{1- t_\al q_\al^{j} q^{(\rho_k-b_-,\al)}}
{1- q_\al^{j}q^{(\rho_k-b_-,\al)}}
\Bigr)
E_{b_-}(\La)\ =\ \si_\ast(q^{b_--\rho_k})\,E_{b_-}(\La)
\notag\\
&= \sum_{w\in W}
w_\La\Bigl(\si_\ast(\La)\,\La_{b^\iota_-}^{-1}\,
\Xi^{(w)}\bigl(q^{b_--\rho_k},\La\bigr)\Bigr)\!
\prod_{\al>0>w(\al)}
\frac{1\!-\!t_\al \La^{-1}_{w(\al)}}
{1 - \La^{-1}_{w(\al)}}\,.
\end{align}

It is instructional to plug in here $\La=q^{b_\#}$
for $b=u_b^{-1}(b_-)$. As we know,
only the term with $w=u_b^{-1}$ in the last sum will
not vanish and $\si_\ast(q^{b_--\rho_k})=
\si_\ast\bigl(q^{u_b(b_\#)}\bigr)=
u_b^{-1}\bigl(\si_\ast(q^{b_\#})\bigr)$; recall that
$b_\#=u_b^{-1}(b_--\rho_k)$.
Then
$u_b^{-1}\bigl(E_b(q^{b_\#})\bigr)=E_b(q^{b_--\rho_k})$,
$\La^{-1}_{w(\al)}=q^{(\rho_k-b_-\,,\,\al)}$ and
we arrive at 
\begin{align*}
E_{b_-}(q^{b_\#})
=\prod_{\al>0>u_b^{-1}(\al)}
\frac{1\!-\!t_\al q^{(\rho_k-b_-\,,\,\al)}}
{1 - q^{(\rho_k-b_-\,,\,\al)}}\,E_b(q^{b_--\rho_k}).
\end{align*}
This is a particular case of the Duality Theorem
from \cite{C4}: 
\begin{align*}
&E_{b}(q^{c_\#})E_{c}(q^{-\rho_k})\ =\ 
E_{c}(q^{b_\#})E_{b}(q^{-\rho_k}),\ b,c\in P.
\end{align*}
Indeed, 
\begin{align*}
&E_{b_-}(q^{-\rho_k})/ 
E_{b}(q^{-\rho_k})=\prod_{\al>0>u_b^{-1}(\al)}
\frac{1\!-\!t_\al q^{(\rho_k-b_-\,,\,\al)}}
{1 - q^{(\rho_k-b_-\,,\,\al)}},
\end{align*}
which is direct from (\ref{macde-eval}).
\smallskip

Note that (\ref{SphPsixx}) becomes very explicit
for $A_1$, which we will omit. Generally, the relation to 
the case of $A_1$ is as follows:
\begin{align*}
\Xi^{(\hbox{\tiny id})}=\tilde{\Xi}_-, \ 
\Xi^{(s)}=\tilde{\Xi}_+\,; \hbox{\, see (\ref{tildeXi})\,}.
\end{align*}
Let us provide two
examples of the coefficients of $\Xi^{w}$
beyond $A_1$. We use the formula from \cite{HHL} with
some help of {\em SAGE}. For $A_2$, let us calculate
the constant term of $\Xi^{(\hbox{\tiny id})}$, which 
is the coefficient of $X_{b_+}$ in the polynomial 
$E_{b_-}\, (b_-\in P_-)$ upon the substitution
$\La=q^{b_--\rho_k}$ for any $b_-\in P_-$. It equals
\begin{align*}
\Xi^{(\hbox{\tiny id})}(X_\al\!\to\!\infty,\al\!>\!0)=
\frac{(1\!-\!t)\bigl(1\!-\!t\!+\!t^2\!-\!t
(\La_{\al_1}^{-1}\!+\!
\La_{\al_2}^{-1}\!-\La_{\al_1\!+\al_2}^{-1})\bigr)}
{(1-\La_{\al_1}^{-1})(1-\La_{\al_2}^{-1})
(1-\La_{\al_1\!+\al_2}^{-1})}.
\end{align*}
The coefficient of  $\Xi^{(\hbox{\tiny id})}$ of $X_{\al_1}^{-1}$
equals
\begin{align*}
\frac{q(1-t)^2}{t(1-q)}\frac{(1-t\La_{\al_2})}{(1-q\La_{\al_2})}
\frac{\bigl(1\!-\!t\!+\!q t^2\!-\!t
(q\La_{\al_1}^{-1}\!+\!
\La_{\al_2}^{-1}\!-\!\La_{\al_1\!+\al_2}^{-1})\bigr)}
{(1-\La_{\al_1}^{-1})(1-\La_{\al_2}^{-1})
(1-\La_{\al_1\!+\al_2}^{-1})}.
\end{align*}
Accordingly, the coefficient  of
$\Xi^{(\hbox{\tiny id})}$ of $X_{\al_2}^{-1}$ is obtained
by the transposition $\al_1\leftrightarrow \al_2$ in the
formula above.
\smallskip

{\sf Symmetrization.}
The passage from the decomposition (\ref{SphPsi})
to that of $F(X,\La)$ can be now performed by means of
the connection formula (\ref{Shinqtsym}). Then we
will arrive at Theorem 4.6 from \cite{Sto2} (in the reduced 
case). We do not provide here the formula for the
radius of convergence of $\Xi^{(w)}(X,\La)$ in terms of $X$. 
Similar to \cite{Sto2}, this can be generally achieved using
the $X\leftrightarrow \La$\~duality and calculating the
distance to the first $\La$\~singular point in the 
$\si$\~decomposition.

The $\si$\~decomposition of $F(X,\La)$ in the case of $t=0$
is directly related to \cite{GiL,BF}  (the 
$q$\~Whittaker case) and the so-called $K$\~theoretic
$J$\~function of flag varieties.
For arbitrary $t$, the connection with the
Laumon spaces was established for $A_n$ in \cite{BFS}; 
see Conjecture 1.8 and especially Proposition 5.11 there.

No geometric interpretation of the $E$\~polynomials
is known for general parameters at the moment. 
It does exist in the following cases. First of all,
$E$\~polynomials coincide with the $p$\~adic spherical
(nonsymmetric) Matsumoto functions for $q=0$. Then
there is a direct link to the level-one Demazure characters
\cite{San,Ion} in the twisted case for $t=0$. More recently,
the global {\em nonsymmetric\,} $q$\~Whittaker function 
and the $E$\~polynomials at $t=\infty$ appeared connected with 
so-called PBW-filtration; see \cite{ChO2,ChF}. 

Generally, 
we have a powerful nonsymmetric machinery of intertwining 
operators \cite{C1,HHL,RY,OS}, which does not exist in the 
symmetric theory, but an obvious lack of geometric 
understanding of the $E$\~polynomials (so far).

\renewcommand\refname{\sc{References}}
\bibliographystyle{unsrt}

\end{document}